\documentclass[12pt, leqno]{article}

\usepackage{amsthm,amsmath}
\usepackage{txfonts}
\usepackage[english]{babel}

\theoremstyle{plain}
 \newtheorem{thm}{Theorem}[section]
 \newtheorem{cor}[thm]{Corollary}
 \newtheorem{lem}[thm]{Lemma}
 \newtheorem{prop}[thm]{Proposition}
 \newtheorem{eg}[thm]{Example}
 \theoremstyle{definition}
 \newtheorem{defn}[thm]{Definition}
 \theoremstyle{remark}
 \newtheorem{rem}[thm]{Remark}

\begin{document}

\title{On Popoviciu Type Formulas for Generalized Restricted Partition Function}

\author{\large Nan LI  and  Sheng CHEN$^{*}$ }

\maketitle

\begin{center}
{Department of Mathematics,   Harbin Institute of Technology,   \\
Harbin 150001,   P.R.China}
\end{center}

\footnotetext{\hspace{-0.5cm} $^*$\hspace*{1mm} Corresponding author.    
\hspace*{3mm}
 {\it E-mail address:} {schen@hit.edu.cn\\ Project 10526016  Supported by National Natural Science Foundation
of China.\\
Project HITC200701 Supported by Science Research Foundation in
Harbin Institute of Technology.}}

\begin{center}
{\bf Abstract\\[1ex]}
\end{center}
Suppose that $a_1(n),a_2(n),\cdots,a_s(n),m(n)$ are integer-valued
polynomials in $n$ with positive leading coefficients. This paper
presents  Popoviciu type formulas for the generalized restricted
partition function $$p_{A(n)}(m(n)):=\#\{(x_1,\cdots,x_s)\in
\mathbb{Z}^{s}:\,\, \,all\,\,x_j\geqslant
0,\,\,x_1a_1(n)+\cdots+x_sa_s(n)=m(n) \}$$ when $s=2$ or $3$. In
either case,
 the formula
implies that the function is an integer-valued quasi-polynomial. The
main result
is proved by  a reciprocity law for a class of fractional part sums and  the theory of generalized
Euclidean division.\\

\noindent 2000 AMS Classification:  Primary  11D45,  Secondary 05A15, 11P99\\
 Key words: \hspace{1mm}
generalized restricted partition function; integer-valued
quasi-polynomial;\par $\quad\quad\quad\,\,\,\,$ reciprocity law;
generalized Euclidean division

\section{Introduction}
We study the generalized restricted partition function
$$p_{A(n)}(m(n)):=\#\{(x_1,\cdots,x_s)\in \mathbb{Z}^{s}:\,\, \,all\,\,x_j\geqslant 0,\,\,x_1a_1(n)+\cdots+x_sa_s(n)=m(n) \}$$
where $a_1(n),a_2(n),\cdots,a_s(n),m(n)$ are integer-valued
polynomials in $n$ with positive leading coefficients, and
$A(n)=\{a_1(n),a_2(n),\cdots,a_s(n)\}$. There are two related
problems (see $\cite{f-d}$, $\cite{coin}$ and $\cite{stanley}$): the
problem of counting the number of lattice points in integer dilates
of the rational polytope and  the linear Diophantine problem of
Frobenius. When $s=2$ and $A(n)\subseteq \mathbb{Z}$, we have the
famous Popoviciu formula (see $\cite{2}$). In the case that $s> 2$,
even when $A(n)=\{a_1(n),a_2(n),\cdots,a_s(n)\}\subseteq
\mathbb{Z}$, explicit formulas for $p_{\{A(n)\}}(m(n))$ have proved
elusive (see $\cite{f-d}$). In $\cite{advance}$, Beck gave a method
to compute refined upper bounds for Frobenious problem when $s=3$.

This paper  aims  to  prove a Popoviciu type formula (see  Theorem
$\ref{p}$ ) for $p_{\{A(n)\}}(m(n))$ in the case that $s=3$, which
implies that it is an integer-valued quasi-polynomial (see Corollary
$\ref{quasi}$). The main result is proved by a reciprocity law (see
($\ref{reci}$) in Theorem $\ref{reciprocity}$) for the following
fractional part sum
$$S(r_1,r_2,r_3;r_4,r_5)=
\sum^{r_2}_{x=r_1}\bigg\{\frac{r_3+r_4x}{r_5}\bigg\}r_5$$ and the
theory of generalized Euclidean division (see Section 3).

\section{A Reciprocity Law for a Class of Fractional Part Sums} Define a class of  fractional part
sums as follows:
\begin{equation}\label{S}
    S(r_1,r_2,r_3;r_4,r_5)=
\sum^{r_2}_{x=r_1}\bigg\{\frac{r_3+r_4x}{r_5}\bigg\}r_5
\end{equation}
 In this section we suppose that $r_1,r_2,r_3,r_4,r_5\in\mathbb{Z}$ and  prove a
reciprocity law and a formula  for $S(r_1,r_2,r_3;r_4,r_5)$ (see
($\ref{reci}$) in Theorem $\ref{reciprocity}$ and ($\ref{S_0}$) in
Corollary $\ref{formula}$).

\begin{lem}\label{4.0.13}
Let $r_1,r_2,r_3,r_4,r_5 \in \mathbb{Z}$. Then
\begin{eqnarray*}
&&\sum_{x=r_1}^{r_2}
\Big\{\frac{r_3+r_4x}{r_5}\Big\}r_5\\
&=&-qr_5+
\sum_{x=r_1}^{\big[\frac{(i_0+1)r_5-r_3}{g_1}\big]}\bigg(r_3+g_1x-i_0r_5\bigg)
+\sum_{i=i_0+1}^{i_t-1}\sum_{x=l_{1i}+1}^{l_{2i}}\bigg(r_3+g_1x-ir_5\bigg)\\
&&+\sum_{x=\big[\frac{i_tr_5-r_3}{g_1}\big]+1}^{r_2}\bigg(r_3+g_1x-i_tr_5\bigg)
\end{eqnarray*}where
\begin{eqnarray*}
      g_1&=& \Big\{\frac{r_4}{r_5}\Big\}r_5,\,\,\,\,i_0 = \bigg[\frac{r_3+g_1r_1}{r_5}\bigg],\,\,\,\, i_t= \bigg[\frac{r_3+g_1r_2}{r_5}\bigg] \\\\
l_{1i} &=& \bigg[\frac{ir_5-r_3}{g_1}\bigg],\,\,\,\,\, \,\,\,l_{2i} =\bigg[\frac{(i+1)r_5-r_3}{g_1}\bigg]\\
 d(n)&=&gcd(r_5,g_1),\,
\widetilde{r_5}=\frac{r_5}{d(n)},\,\widetilde{g_1}=\frac{g_1}{d(n)},\,\widetilde
{r_3}=\frac{r_3}{d(n)}\\
q&=& \left\{
  \begin{array}{ll}
    \big[\frac{i_t-i_0-1-r}{\widetilde{g_1}}\big]+1, & \hbox{$d(n)\mid r_3$;} \\
    0, & \hbox{$d(n)\nmid r_3$.}
  \end{array}
\right.\\
r&=&\bigg\{\frac{(\widetilde{r_3}-(i_0+1)\widetilde{r_5})\widetilde{r_5}^{-1}}{\widetilde{g_1}}\bigg\}\widetilde{g_1}
\end{eqnarray*}
in which $\widetilde{r_5}^{-1}\in \mathbb{Z}$ and
$\widetilde{r_5}^{-1}\widetilde{r_5} \equiv 1$ $(mod\,\,
\widetilde{g_1})$.
\end{lem}
\begin{proof}
From the definition of $g_1$, we have $r_4\equiv g_1(mod \,r_5)$,
and therefore
$$\Big\{\frac{r_3+r_4x}{r_5}\Big\}r_5=\Big\{\frac{r_3+g_1x}{r_5}\Big\}r_5$$
By  Euclidean division,  we have
  $$\bigg\{\frac{r_3+g_1r_1}{r_5}\bigg\}r_5 =
  (r_3+g_1r_1)-r_5i_0,\,\,\,where\,\,i_0=\bigg[\frac{r_3+g_1r_1}{r_5}\bigg]$$
  and
 $$\bigg\{\frac{r_3+g_1r_2}{r_5}\bigg\}r_5 =
 (r_3+g_1r_2)-r_5i_t,\,\,\,where\,\,i_t=\bigg[\frac{r_3+g_1r_2}{r_5}\bigg]$$
 For arbitrary integer $x$ with $r_1< x< r_2$, it is easy to check
that
$$\bigg\{\frac{r_3+g_1x}{r_5}\bigg\}r_5=r_3+g_1x-r_5i$$
if and only if
$$\frac{r_5i-r_3}{g_1}\leqslant x < \frac{r_5(i+1)-r_3}{g_1}$$

(I) First we consider the case that $gcd(r_5,g_1)\nmid r_3$, which
implies that $\frac{r_5i-r_3}{g_1}$ can not be an integer for any
$i\in\mathbb{Z}$. Then we have
$$\bigg\{\frac{r_3+g_1x}{r_5}\bigg\}r_5=r_3+g_1x-r_5i$$
if and only if
$$l_{1i}+1=\bigg[\frac{r_5i-r_3}{g_1}\bigg]+1\leqslant x \leqslant \bigg[ \frac{r_5(i+1)-r_3}{g_1}\bigg]=l_{2i}$$
Let $S=\sum_{x=r_1}^{r_2} \Big\{\frac{r_3+r_4x}{r_5}\Big\}r_5$. Then
\begin{equation}\label{*}
    S=\sum_{x=r_1}^{\big[\frac{(i_0+1)r_5-r_3}{g_1}\big]}\bigg(r_3+g_1x-i_0r_5\bigg)
+\sum_{i=i_0+1}^{i_t-1}\sum_{x=l_{1i}+1}^{l_{2i}}\bigg(r_3+g_1x-ir_5\bigg)
+\sum_{x=\big[\frac{i_tr_5-r_3}{g_1}\big]+1}^{r_2}\bigg(r_3+g_1x-i_tr_5\bigg)
\end{equation}

Since $gcd(r_5,g_1)\nmid r_3$, for any $i_0+1\leqslant i\leqslant
i_t$, we have
\begin{equation}\label{i-1}
\bigg\{\frac{r_3+g_1x}{r_5}\bigg\}r_5=r_3+g_1x-r_5(i-1),\,\,\,\,if\,\,\,\,x=\bigg[\frac{r_5i-r_3}{g_1}\bigg]
\end{equation}

(II) Now we consider the case  that $gcd(r_5,g_1)\mid r_3$. Let
$\overline{i}\in \mathbb{Z}$ such that
$\frac{r_5\overline{i}-r_3}{g_1}$ is an integer. Then we have
\begin{equation}\label{i}
\bigg\{\frac{r_3+g_1x}{r_5}\bigg\}r_5=r_3+g_1x-r_5\overline{i},\,\,\,\,if\,\,\,\,x=\bigg[\frac{r_5\overline{i}-r_3}{g_1}\bigg]
\end{equation}
Note that in this case,
$$i=\overline{i},\,\,\,k=\frac{r_5\overline{i}-r_3}{g_1}\in \mathbb{Z}$$ satisfy the equation $$r_5i-g_1k=r_3$$Comparing
($\ref{i-1}$) and ($\ref{i}$), we know that in order to get the
formula in this case, we  need only to add a term $-qr_5$ to
($\ref{*}$), where $q$ is the number of integer solutions to
$r_5i-g_1k=r_3$ with $i_0+1 \leqslant i\leqslant i_t$.

Since $q$ can be calculated easily, we have the formula in this
lemma.
\end{proof}

From the lemma above, the following reciprocity law for the
fractional part sum $S(r_1,r_2,r_3;r_4,r_5)$ follows.
\begin{thm}\label{reciprocity}Let $r_1,r_2,r_3,r_4,r_5\in \mathbb{Z}$. Then
\begin{equation}\label{reci}
S(r_1,r_2,r_3;r_4,r_5)+\frac{r_5}{g_1}S(r_1^{'},r_2^{'},-r_3;r_5,g_1)=A
\end{equation}
where $$g_1=\bigg\{\frac{r_4}{r_5}\bigg\}r_5$$
$$r_1^{'}=
               \bigg[\frac{r_3+g_1r_1}{r_5}\bigg]+2,\,\,\,r_2^{'} =
               \bigg[\frac{r_3+g_1r_2}{r_5}\bigg]-1$$
\begin{eqnarray*}
  A &=&
  (r_2-r_1+1)r_3+\frac{g_1}{2}(r_2+r_1)(r_2-r_1+1)\\
 &&-r_5\Bigg(-\bigg[\frac{(r_1^{'}-1)r_5-r_3}{g_1}\bigg]+(r_1^{'}-2)(-r_1+1)\\
 &&
 +q-\bigg[\frac{(r_2^{'}+1)r_5-r_3}{g_1}\bigg]+(r_2^{'}+1)r_2\\
  &&+\frac{r_3}{g_1}(r_2^{'}-r_1^{'}+1)-\frac{r_5}{2g_1}(r_2^{'}+r_1^{'})(r_2^{'}-r_1^{'}+1)\Bigg).\\
  q&=& \left\{
  \begin{array}{ll}
    0, & \hbox{$gcd(r_5,g_1)\nmid r_3$;} \\
    \bigg[\frac{r_2^{'}-r_1^{'}+2-r}{\frac{g_1}{gcd(r_5,g_1)}}\bigg]+1, & \hbox{$gcd(r_5,g_1)\mid r_3$.}
  \end{array}
\right.
\\r&=&\bigg\{\frac{\big(\frac{r_3}{gcd(r_5,g_1)}-(r_1^{'}-1)\frac{r_5}{gcd(r_5,g_1)}\big)\big(\frac{r_5}{gcd(r_5,g_1)}\big)^{-1}}
{\frac{g_{1}}{gcd(r_5,g_1)}}\bigg\}\frac{g_{1}}{gcd(r_5,g_1)}
\end{eqnarray*}
in which $\big(\frac{r_5}{gcd(r_5,g_1)}\big)^{-1}\in\mathbb{Z}$ and
$\big(\frac{r_5}{gcd(r_5,g_1)}\big)^{-1}\frac{r_5}{gcd(r_5,g_1)}\equiv
1$ $(mod\,\,\frac{g_{1}}{gcd(r_5,g_1)})$.
\end{thm}

Now we can  give an explicit formula to compute
$S(r_1,r_2,r_3;r_4,r_5)$ using the reciprocity law. Put
$$g_{-1}=r_4,\,g_0=r_5$$and $$g_k=\bigg\{\frac{g_{k-2}}{g_{k-1}}\bigg\}g_{k-1}$$ for $1\leqslant k \leqslant
h$. Suppose that $g_h\mid g_{h-1}$. For simplicity, we write
$$g_{-1}\rightarrow g_0\rightarrow g_1\rightarrow g_2\rightarrow
\cdots \rightarrow g_{h-1}\rightarrow g_h\rightarrow 0$$
 For $0\leqslant k \leqslant h$, put
$$S_k =S(r_1^{(k)},r_2^{(k)},(-1)^{k}r_3;g_{k-1},g_{k})$$
where $$r_1^{(0)}=r_1,\,\, r_2^{(0)}=r_2$$
       and \begin{eqnarray*}
r_1^{(k)} &=&
               \bigg[\frac{(-1)^{k-1}r_3+g_kr_1^{(k-1)}}{g_{k-1}}\bigg]+2\\
                r_2^{(k)} &=&
               \bigg[\frac{(-1)^{k-1}r_3+g_kr_2^{(k-1)}}{g_{k-1}}\bigg]-1
             \end{eqnarray*}
for $1\leqslant k \leqslant h$. Then we have the following formula
for $S(r_1,r_2,r_3;r_4,r_5)$ using  the reciprocity law in Theorem
$\ref{reciprocity}$ successively.

\begin{cor}\label{formula}
Let
$$S_0=S(r_1,r_2,r_3;g_{-1},g_0)=\sum^{r_2}_{x=r_1}\bigg\{\frac{r_3+g_{-1}x}{g_0}\bigg\}g_0$$
where $r_1,r_2,r_3,r_4,r_5\in \mathbb{Z}$. Using the notations
above, we have
 \begin{equation}\label{S_0}
S_0=A_0+g_0\sum^{h-1}_{i=1}(-1)^{i}\frac{A_{i}}{g_i}+(-1)^{h}\frac{g_{0}}{g_h}S_h
 \end{equation}
where
 $S_h$ can be
calculated directly as follows (since $g_h\mid g_{h-1}$)
$$S_h=\sum^{r_2^{(h)}}_{x=r_1^{(h)}}\bigg\{\frac{(-1)^{h}r_3+g_{h-1}x}{g_h}\bigg\}g_h
=\sum^{r_2^{(h)}}_{x=r_1^{(h)}}\bigg\{\frac{(-1)^{h}r_3}{g_h}\bigg\}g_h=(r_2^{(h)}-r_1^{(h)}+1)\bigg\{\frac{(-1)^{h}r_3}{g_h}\bigg\}g_h$$
and
\begin{eqnarray*}
  A_{k-1} &=&
  (r_2^{(k-1)}-r_1^{(k-1)}+1)(-1)^{k-1}r_3+\frac{g_k}{2}(r_2^{(k-1)}+r_1^{(k-1)})(r_2^{(k-1)}-r_1^{(k-1)}+1)\\
 &&-g_{k-1}\Bigg(-\bigg[\frac{(r_1^{(k)}-1)g_{k-1}+(-1)^{k}r_3}{g_k}\bigg]+(r_1^{(k)}-2)(-r_1^{(k-1)}+1)\\
 &&
 +q_{k-1}-\bigg[\frac{(r_2^{(k)}+1)g_{k-1}+(-1)^{k}r_3}{g_k}\bigg]+(r_2^{(k)}+1)r_2^{(k-1)}\\
  &&+\frac{(-1)^{k-1}r_3}{g_k}(r_2^{(k)}-r_1^{(k)}+1)-\frac{g_{k-1}}{2g_k}(r_2^{(k)}+r_1^{(k)})(r_2^{(k)}-r_1^{(k)}+1)\Bigg)\\
  q_{k-1}&=& \left\{
  \begin{array}{ll}
    0, & \hbox{$g_h\nmid r_3$;} \\
    \bigg[\frac{r_2^{(k)}-r_1^{(k)}+2-r}{\frac{g_k}{g_h}}\bigg]+1, & \hbox{$g_h\mid r_3$.}
  \end{array}
\right.\\
r&=&\bigg\{\frac{\big(\frac{(-1)^{k-1}r_3}{g_h}-(r_1^{(k)}-1)\frac{g_{k-1}}{g_h}\big)\big(\frac{g_{k-1}}{g_h}\big)^{-1}}
{\frac{g_{k}}{g_h}}\bigg\}\frac{g_{k}}{g_h}
\end{eqnarray*}
in which  $\big(\frac{g_{k-1}}{g_h}\big)^{-1}\in \mathbb{Z}$ and
$\big(\frac{g_{k-1}}{g_h}\big)^{-1}\frac{g_{k-1}}{g_h}\equiv 1$
$(mod\,\,\frac{g_{k}}{g_h})$.
\end{cor}
\begin{proof}
By the reciprocity law  in Theorem $\ref{reciprocity}$ for the
fractional part sums
$$S_{k-1}=S(r_1^{(k-1)},r_2^{(k-1)},(-1)^{k-1}r_3;g_{k-2},g_{k-1})$$
we obtain
$$S(r_1^{(k-1)},r_2^{(k-1)},(-1)^{k-1}r_3;g_{k-2},g_{k-1})+\frac{g_{k-1}}{g_k}S(r_1^{(k)},r_2^{(k)},(-1)^{k}r_3;g_{k-1},g_{k})=A_{k-1}$$
Then this corollary follows.
\end{proof}
\section{The Theory of Generalized Euclidean Division and GCD}
To prove the main result in this paper, we need the theory of
Euclidean division and GCD    for the ring $R$ of integer-valued
quasi-polynomials (see Definition $\ref{2.0.1}$ and Proposition
$\ref{R}$). Such a theory has been discussed in $\cite{mine}$. In
this section we list related definitions and results without proofs.
It should be pointed out that the definition of integer-valued
quasi-polynomials defined here is slightly different from that in
$\cite{mine}$, i.e., integer-valued quasi-polynomials in Definition
$\ref{2.0.1}$ have lower boundaries while in $\cite{mine}$ they do
not. However, the following results and proofs hold in spite of this
modification.

\begin{defn}[Integer-valued quasi-polynomial, similar to Definition 3 in $\cite{mine}$]\label{2.0.1}
We call function $f:\mathbb{N}\to\mathbb{Z}$ an integer-valued
quasi-polynomial, if there exist positive integers $T,C$, and
polynomials $f_i(x)\in \mathbb{Z}[x]$ ($i=0,1,\cdots,T-1$), such
that if $n>C$ and $n=Tm+i$ ($m\in \mathbb{N}$), we have
$f(n)=f_i(m)$. We call $(T, C, \{f_i(x)\}_{i=0}^{T-1})$ a
\emph{representation} of $f(x)$ and write
$$f(x)=(T, C, \{f_i(x)\}_{i=0}^{T-1})$$ Then $\max\{degree(f_i(x))|
i=0,  1,  \cdots,  T-1\}$, $T$ and $C$ are called the  degree,
period and lower boundary of this representation respectively.
\end{defn}

\begin{prop}[Definition of the ring $R$, see Proposition 6 in $\cite{mine}$]\label{R}The set of all integer-valued quasi-polynomials,   denoted
by $R$,   with pointwisely defined addition and multiplication, is a
commutative ring with identity.
\end{prop}

\begin{defn}[See Definition 9 in $\cite{mine}$]\label{2. 0. 2}
               Let $r(x)\in R$ and $r(x)=(T,C,
 \{r_i(x)\}_{i=0}^{T-1})$. We shall say $r(x)$ is nonnegative and  write $r(x)\succcurlyeq  0$,   if it
satisfies the following equivalent conditions:
\begin{enumerate}
\item[(a)]
for every $i=0,  1,  \cdots,  T-1$, $r_i(x)=0$ or its leading
coefficient is positive;
          \item[(b)] there exists $C_1\in \mathbb{Z}$,   such that for every integer $n>C_1$,  we have $r(n) \geqslant 0$.
          \end{enumerate}
We shall say $r(x)$ is strictly positive and  write $r(x)\succ 0$,
if $r(x)=(T,C,
 \{r_i(x)\}_{i=0}^{T-1})$ satisfies the following  condition:

 $(a^{'})$ for every $i=0,  1,  \cdots, T-1$,
the leading coefficient of $r_i(x)$ is positive.\\
 {We write $
f(x)\preccurlyeq g(x)$  if $ g(x)-  f(x)\succcurlyeq 0$}.
\end{defn}
\begin{defn}[See Definition 11 in $\cite{mine}$]\label{2. 0. 20}
Let $r(x) \in \mathbb{Z}[x] $,   define a function $| \cdot  |:$
$\mathbb{Z}[x]\rightarrow \mathbb{Z}[x]$ as follows:
$$|r(x)|=\left\{
          \begin{array}{ll}
            r(x),   & if\,  \,  \hbox{$r(x)\succ 0$} \\
            -r(x),   & if \,  \,   \hbox{$r(x)\prec 0$} \\
            0,   & if \,  \,  \hbox{$r(x)=0$}
          \end{array}
        \right. $$
\end{defn}

\begin{thm}[Generalized Euclidean division, see Theorem 12 in $\cite{mine}$]\label{2. 0. 3}
Let $f(x),  g(x)\in \mathbb{Z}[x]$ and $g(x)\neq 0$.  Then there
exist unique $P(x),  r(x)\in R$ such that
$$f(x)=P(x)g(x)+r(x)\quad\,  where\quad
0\preccurlyeq r(x)\prec| g(x)|$$ In this situation,  we write
$P(x)=quo(f(x), g(x))$ and $r(x)=rem(f(x),g(x))$. We call  the
inequality $0\leqslant r_1(x)<|g(x)|$ the remainder condition of
this division.\end{thm}
\begin{rem}[Similar to Remark 13 in $\cite{mine}$]\label{2. 0. 01}This division almost coincides with
the division in $\mathbb{Z}$ pointwisely in the following sense.
Let $f(x),  g(x)\in \mathbb{Z}[x]$ and
$$f(x)=P(x)g(x)+r(x),\quad\,  where\quad
0\preccurlyeq r(x)\prec| g(x)|$$By Definition $\ref{2. 0. 2}$,   the
inequality $0\leqslant r(x)<|g(x)|$ will give an integer $C$  such
that for all $n
>C$,
$$\big[\frac{f(n)}{g(n)}\big]=P(n)\,  ,  \,  \big\{\frac{f(n)}{g(n)}\big\}g(n)=r(n)$$
If $r(x)=0$, we have
$$r(n)=\big\{\frac{f(n)}{g(n)}\big\}g(n)=0$$ for every $n\in \mathbb{N}$. The opposite also holds,
i.e. if  $r(n)=0$ for every $n\in \mathbb{N}$, then $r(x)=0$ as an
integer-valued quasi-polynomial. So $ rem(f(x),  g(x))=0$ if and
only if for every $n\in \mathbb{N}$,  $ rem(f(n),  g(n))=0$.
\end{rem}

\begin{eg}[Similar to Example 14 in $\cite{mine}$]\label{eg1}The following is  an example to
illustrate the relation between Euclidean division in
$\mathbb{Z}[x]$ and on $\mathbb{Z}$.

 When $n>1$,
$$\bigg[\frac{n^{2}}{2n+1}\bigg]=\left\{
                             \begin{array}{ll}
                               m-1,   & \hbox{$n=2m$;} \\
                               m-1,   & \hbox{$n=2m-1$. }
                             \end{array}
                           \right. $$
$$\bigg\{\frac{n^{2}}{2n+1}\bigg\}(2n+1)=\left\{
                             \begin{array}{ll}
                               3m+1,   & \hbox{$n=2m$;} \\
                               m,   & \hbox{$n=2m-1$. }
                             \end{array}
                           \right.
$$
\end{eg}
 In order to use successive division in $R$,  we define generalized Euclidean
Division for $f(x),g(x)\in R$ as follows, similar to  $\cite{mine}$.
Suppose that $T_0$ is the least common period of $f(x)$, $g(x)$,
such that $f(x)=(T_0, C_1,\{f_i(x)\}_{i=0}^{T_0-1})$ and
$g(x)=(T_0,C_2, \{g_i(x)\}_{i=0}^{T_0-1})$.  Based on generalized
division in $\mathbb{Z}[x]$ (see Theorem $\ref{2. 0. 3}$),   we can
define $quo(f(x),  g(x))$ and $rem(f(x), g(x))$ as follows (denoted
by $P(x)$ and $r(x)$ respectively): when $n=Tm+i$
$$P(n)=\left\{
         \begin{array}{ll}
           quo(f_i(m),  g_i(m)),   &  \hbox{ $g_i(m)\neq 0$} \\
           0,   &  \hbox{$g_i(m)=0$}
         \end{array}
       \right.,
r(n)=\left\{
         \begin{array}{ll}
           rem(f_i(m),  g_i(m)),   & \hbox{ $g_i(m)\neq 0$} \\
           f_i(m),   & \hbox{$g_i(m)=0$}
         \end{array}
       \right. $$
Then it is easy to check that $P(x),   r(x)\in R$ and
\begin{equation}\label{fg}
f(x)=R(x)g(x)+r(x)
\end{equation}
This will be called  the generalized Euclidean algorithm on the ring
of integer-valued quasi-polynomials.

By successive division in $R$, we can develop generalized GCD theory
(see Definition $\ref{2. 0. 006}$), similar to the case of
$\mathbb{Z}$ (see \cite{Bhu}).

\begin{defn}[Divisor, similar to Definition 15 in $\cite{mine}$]\label{divisor}Suppose that $f(x),  g(x)\in R$ and for every $n\in \mathbb{N}$,
$g(n)\neq 0$.  Then by Remark $\ref{2. 0. 01}$,   the following two
conditions are equivalent:
\begin{enumerate}
  \item[(1)]rem($f(x),  g(x)$)=0;
  \item[(2)] for every $x\in \mathbb{N}$,   $g(x)$ is a divisor of $f(x)$.
\end{enumerate}
If  the two conditions are satisfied,   we shall call $g(x)$ a
divisor of $f(x)$ and write $g(x)\mid f(x)$.
\end{defn}
\begin{defn}[Quasi-rational function]\label{quasi-rational}
A  function $f:\mathbb{N}\to\mathbb{Q}$  is quasi-rational, if there
exist positive integers $T,C$, and rational functions
$\frac{f_i(x)}{g_i(x)}\in \mathbb{Q}(x)$, where $f_i(x),g_i(x)\in
\mathbb{Z}[x]$ ($i=0,1,\cdots,T-1$), such that if $n>C$ and $n=Tm+i$
($m\in \mathbb{N}$), we have $f(n)=\frac{f_i(m)}{g_i(m)}$.
\end{defn}
By the equivalence of the  two conditions in Definition
$\ref{divisor}$, we have the following property for quasi-rational
functions.

\begin{prop}\label{rational}
Let $f(x)$ be a  quasi-rational function. If for every $n\in
\mathbb{N}$, $f(n)\in \mathbb{Z}$, then $f(x)$ is an integer-valued
quasi-polynomial.
\end{prop}
Besides, by the equivalence of the two conditions in Definition
$\ref{divisor}$,  we have the following proposition,  similar to the
situation in $\mathbb{Z}$.
 \begin{prop}[See Proposition 16 in $\cite{mine}$]\label{2. 0. 6}Let
$g(x),   f(x)\in $R.  If $f(x)\mid g(x)$ and $g(x)\mid f(x)$,   we
have $f(x)=\varepsilon g(x)$,   where $\varepsilon $ is an
invertible element in $R$ (see Proposition 7 in $\cite{mine}$ for
description of invertible elements in $R$).
\end{prop}

\begin{defn} [Generalized GCD, see Definition 17 in $\cite{mine}$]\label{2. 0. 006}Let $f_1(x),  f_2(x),  \cdots,  f_s(x),  d(x)\\\in
R$.
\begin{enumerate}
  \item[(1)] We call $d(x)$  a
common divisor of $f_1(x),  f_2(x),  \cdots,  f_s(x)$,   if we have
$d(x)\mid f_k(x)$ for every $k=1,  2,  \cdots,  s$.
  \item[(2)]We call
$d(x)$  a greatest common divisor  of $f_1(x),  f_2(x),  \cdots,
f_s(x)$ if $d(x)$ is a common divisor  of $f_1(x),  f_2(x),  \cdots,
f_s(x)$ and  for any common divisor $p(x)\in R$ of $f_1(x), f_2(x),
\cdots, f_s(x)$,  we have $p(x)\mid d(x)$.
\end{enumerate}
\end{defn}
\begin{rem}[See Remark 18 in $\cite{mine}$]\label{2. 0. 014} Suppose that both $d_1(x)$ and   $ d_2(x)$  are
 greatest common divisors of $f_1(x),  f_2(x), \cdots,  f_s(x)$.
  Then we have $d_1(x)\mid d_2(x)$ and $d_2(x)\mid d_1(x)$. Thus, by
Proposition $\ref{2. 0. 6}$,   we have $d_1(x)=\varepsilon d_2(x)$,
where $\varepsilon $ is an invertible element in $R$.  So we have a
unique GCD $d(x)\in R$ for $f_1(x), f_2(x),  \cdots,  f_s(x)$ such
that $d(x)\succcurlyeq  0$ and write it as $ggcd(f_1(x), f_2(x),
\cdots, f_s(x))$.
\end{rem}
\begin{thm}[See Theorem 21 in $\cite{mine}$]\label{2. 0. 0099}Let $f_1(x),  f_2(x),  \cdots,  f_s(x)\in R$.
Then there exist $d(x),  u_i(x)\in R$ ($ i=1,  2,  \cdots,  s$),
such that $d(x)=ggcd(f_1(x),  f_2(x),  \cdots,  f_s(x))$ and
$$f_1(x)u_1(x)+f_2(x)u_2(x)+\cdots+f_s(x)u_s(x)=d(x)$$
\end{thm}

\begin{defn}\label{3.0.8}Let $a_1(x),a_2(x)$ be integer-valued quasi-polynomials in $x$.
Suppose that $ggcd(a(x),a_2(x))=1$. By Theorem $\ref{2. 0. 0099}$,
there exist two integer-valued quasi-polynomials $u_1(x),u_2(x)$,
such that $a_1(x)u_1(x)+a_2(x)u_2(x)=1$. Then we shall call $u_1(x)$
an \emph{inverse} of $a_1(x)$ mod $a_2(x)$ and $u_2(x)$ an
\emph{inverse} of $a_2(x)$ mod  $a_1(x)$, denoted by $a_1(x)^{-1}$
and $a_2(x)^{-1}$ respectively.
\end{defn}

Now we get a reciprocity law and a formula similar to the results in
Section 2.

\begin{rem}\label{general}Let $r_1,r_2,r_3,r_4,r_5$  be integer-valued quasi-polynomials instead of
integers in $S(r_1,r_2,r_3;r_4,r_5)$ (see ($\ref{S}$)). Then Lemma
$\ref{4.0.13}$, Theorem $\ref{reciprocity}$ and Corollary
$\ref{formula}$ still hold if we use generalized Euclidean division
and  GCD for the ring $R$  of integer-valued quasi-polynomials
 instead of instead of the classical ones for $\mathbb{Z}$.
\end{rem}

\section{A Formula for $p_{\{a_1(n),a_2(n),a_3(n)\}}(m(n))$}

First we want to apply Popoviciu formula to
$$p_{\{a_1(n),a_2(n)\}}(m(n)):=\#\{(x_1,x_2)\in \mathbb{Z}^{2}:\,\, \,all\,\,x_j\geqslant 0,\,\,x_1a_1(n)+x_2a_2(n)=m(n) \}$$
where $a_1(n),a_2(n),m(n)$ are integer-valued quasi-polynomials (see
Definition $\ref{2.0.1}$) with positive leading coefficients. Let
$d(n)=ggcd(a_1(n),a_2(n))$ (see Remark $\ref{2. 0. 014}$ for
notation of $ggcd$). Define
$$\chiup(m(n),d(n))=\left\{
  \begin{array}{ll}
    1, & \hbox{$d(n)\mid m(n)$} \\
    0, & \hbox{$d(n)\nmid m(n)$}
  \end{array}
\right.$$ It is obvious that
$$p_{\{a_1(n),a_2(n)\}}(m(n))=\chiup(m(n),d(n))p_{\{\frac{a_1(n)}{d(n)},\frac{a_2(n)}{d(n)}\}}\bigg(\frac{m(n)}{d(n)}\bigg)$$
So it is sufficient to consider the case  that $d(n)=1$. In this
case, by generalized Euclidean division (see Theorem $\ref{2. 0.
3}$), generalized GCD algorithm (see Definition $\ref{2. 0. 006}$
and Theorem $\ref{2. 0. 0099}$) and inverse (see Definition
$\ref{3.0.8}$), we can apply Popoviciu formula to
$p_{\{a_1(n),a_2(n)\}}(m(n))$.
\begin{lem}\label{3.0.9}
If $ggcd(a(n),b(n))=1$, then
$$p_{\{a_1(n),a_2(n)\}}(m(n))=\frac{m(n)}{a_1(n)a_2(n)}-\bigg\{\frac{a_1(n)^{-1}m(n)}{a_2(n)}\bigg\}-
\bigg\{\frac{a_2(n)^{-1}m(n)}{a_1(n)}\bigg\}+1$$ is an
integer-valued quasi-polynomial, where $a_1(n)^{-1}$ is an inverse
of $a_1(n)$ mod $a_2(n)$ and $a_2(n)^{-1}$ is an inverse of $a_2(n)$
mod $a_1(n)$.
\end{lem}
\begin{proof}
By Theorem $\ref{2. 0. 3}$,
$\bigg\{\frac{a_1(n)^{-1}m(n)}{a_2(n)}\bigg\}$ and
$\bigg\{\frac{a_2(n)^{-1}m(n)}{a_1(n)}\bigg\}$ are quasi-rational
(see Definition $\ref{quasi-rational}$). Note that for any $n\in
\mathbb{N}$, $p_{\{a_1(n),a_2(n)\}}(m(n))\in \mathbb{N}$. Therefore,
by Proposition $\ref{rational}$, $p_{\{a_1(n),a_2(n)\}}(m(n))$ is an
integer-valued quasi-polynomial.\end{proof}

Now we turn to
$$p_{\{a_1(n),a_2(n),a_3(n)\}}(m(n)):=\#\{(x_1,x_2,x_3)\in \mathbb{Z}^{3}:\,\, \,all\,\,x_j\geqslant 0,\,\,x_1a_1(n)+x_2a_2(n)+x_3a_3(n)=m(n) \}$$
where $a_1(n),a_2(n),a_3(n),m(n)$ are integer-valued
quasi-polynomials (see Definition $\ref{2.0.1}$) with positive
leading coefficients. To deal with the case of
$a_1(n),a_2(n),a_3(n)$ being not pairwise prime, we need the
following lemma.
\begin{lem}\label{prime}
Let
$$d_0(n)=ggcd(a_1(n),a_2(n),a_{3}(n)),\,\,\,\,d(n)=ggcd(a_1^{'}(n),a_2^{'}(n))$$(see Remark $\ref{2. 0. 014}$  for notation of $ggcd$)
where
$$m^{'}(n)=\frac{m(n)}{d_0(n)},\,\,\,\,a_i^{'}(n)=\frac{a_i(n)}{d_0(n)},
\,\,\,\,i=1,2,3$$ Define
$$\chiup(m(n),d_0(n))=\left\{
  \begin{array}{ll}
    1, & \hbox{$d_0(n)\mid m(n)$;} \\
    0, & \hbox{$d_0(n)\nmid m(n)$.}
  \end{array}
\right.$$ Then we have
$$p_{\{a_1(n),a_2(n),a_3(n)\}}(m(n))=\chiup(m(n),d_0(n))\sum_{k=0}^{\big[\frac{\widetilde{m^{'}(n)}}{a_3^{'}(n)}\big]}
p_{\{\widetilde{a_1^{'}(n)},\widetilde{a_2^{'}(n)}\}}
(\widetilde{m^{'}(n)}-a_3^{'}(n)k)$$ where
$$\widetilde{m^{'}(n)}=\frac{m^{'}(n)-a_3^{'}(n)x_0}{d(n)},\,\,\widetilde{a_i^{'}(n)}=\frac{a_i^{'}(n)}{d(n)}\,\,\,for\,\,i=1,2$$and
$$x_0=\bigg\{\frac{m^{'}(n)a_3^{'}(n)^{-1}}{d(n)}\bigg\}d(n)$$in
which $a_3^{'}(n)^{-1}$ is an inverse of $a_3^{'}(n)$ mod $d(n)$.
\end{lem}
\begin{proof}It is obvious that
$$p_{\{a_1(n),a_2(n),a_3(n)\}}(m(n))=\chiup(m(n),d_0(n))p_{\{a_1^{'}(n),a_2^{'}(n),a_3^{'}(n)\}}(m^{'}(n))$$
and $$ggcd(a_1^{'}(n),a_2^{'}(n),a_3^{'}(n))=1.$$ Since when
$d_0(n)\nmid m(n)$, $p_{\{a_1(n),a_2(n),a_3(n)\}}(m(n))=0$, it
surfies to consider the case that $d_0(n)\mid m(n)$, i.e. $m^{'}(n)$
is an integer-valued quasi-polynomial. Note that
$$\widetilde{a_1^{'}(n)}x_1+\widetilde{a_{2}^{'}(n)}x_{2}=\frac{m^{'}(n)-a_3^{'}(n)x_3}{d(n)}.$$
Denote $p=\frac{m^{'}(n)-a_3^{'}(n)x_3}{d(n)}$, then $x_3$ and $p$
are non-negative integers satisfying the equation
$$a_3^{'}(n)x_3+d(n)p=m^{'}(n).$$ Denote the minimal possible value of $x_3$
by $x_0$ and
$$x_0=\bigg\{\frac{m^{'}(n)a_3^{'}(n)^{-1}}{d(n)}\bigg\}d(n)$$where
$a_3^{'}(n)^{-1}$  is an inverse of $a_3^{'}(n)$ mod $d(n)$. Then we
have $x_3=x_0+kd(n)$ for $k\in\mathbb{Z}$ and $0\leqslant k
\leqslant \big[\frac{m^{'}(n)-a_3^{'}(n)x_0}{a_3^{'}(n)d(n)}\big]$.
Thus
\begin{eqnarray*}
     p_{\{a_1^{'}(n),a_2^{'}(n),a_3^{'}(n)\}}(m^{'}(n))&=&
   \sum_{k=0}^{\big[\frac{m^{'}(n)-a_3^{'}(n)x_0}{a_3^{'}(n)d(n)}\big]}
   p_{\{\widetilde{a_1^{'}(n)},\widetilde{a_2^{'}(n)}\}}\bigg(\frac{m^{'}(n)-a_3^{'}(n)(x_0+k d(n))}{d(n)}\bigg)\\
&=&
   \sum_{k=0}^{\big[\frac{\widetilde{m^{'}(n)}}{a_3^{'}(n)}\big]}
p_{\{\widetilde{a_1^{'}(n)},\widetilde{a_2^{'}(n)}\}}
(\widetilde{m^{'}(n)}-a_3^{'}(n)k)
\end{eqnarray*}
where $\widetilde{m^{'}(n)}=\frac{m^{'}(n)-a_3^{'}(n)x_0}{d(n)}$ is
an integer-valued quasi-polynomial.
\end{proof}
Then by Popoviciu formula (see Lemma $\ref{3.0.9}$) and the formula
for a class of fractional part sums (see Corollary $\ref{formula}$),
we have a Popviciu type formula for the generalized restricted
partition function $p_{\{a_1(n),a_2(n),a_3(n)\}}(m(n))$ as follows.
\begin{thm}\label{p}Let $a_1(n),a_2(n),a_3(n),m(n)$ be integer-valued
quasi-polynomials (see Definition $\ref{2.0.1}$) with positive
leading coefficients. For simplicity, we denote
$\widetilde{m^{'}(n)},\widetilde{a_1^{'}(n)},\widetilde{a_2^{'}(n)},a_3^{'}(n)$
in Lemma $\ref{prime}$ and $
\widetilde{a_1^{'}(n)}^{-1},\widetilde{a_2^{'}(n)}^{-1}$ by
$m,a_1,a_2,a_3$ and $a_1^{-1},a_2^{-1}$ respectively. Then we have
$$p_{\{a_1(n),a_2(n),a_3(n)\}}(m(n))=C-\frac{1}{a_2}S(0,\big[\frac{m}{a_3}\big],ma_1^{-1};-a_3a_1^{-1},a_2)
-\frac{1}{a_1}S(0,\big[\frac{m}{a_3}\big],ma_2^{-1};-a_3a_2^{-1},a_1)$$
where
$$C=\sum^{[\frac{m}{a_3}]}_{x=0}\bigg(\frac{m-a_3x}{a_1a_2}+1\bigg)$$
(Formulas for
$S(0,\big[\frac{m}{a_3}\big],ma_1^{-1};-a_3a_1^{-1},a_2)$ and
$S(0,\big[\frac{m}{a_3}\big],ma_2^{-1};-a_3a_2^{-1},a_1)$ are
referred to  Corollary $\ref{formula}$ and Remark $\ref{general}$).
\end{thm}

\begin{cor}\label{quasi}Let $a_1(n),a_2(n),a_3(n),m(n)$ be integer-valued
quasi-polynomials with positive leading coefficients. Then
$p_{\{a_1(n),a_2(n),a_3(n)\}}(m(n))$ is an integer-valued
quasi-polynomial.
\end{cor}
\begin{proof}The formula for $p_{\{a_1(n),a_2(n),a_3(n)\}}(m(n))$ given in Theorem
$\ref{p}$ is expressed in terms of quotients and remainders of
generalized Euclidean division. The result follows from Theorem
$\ref{2. 0. 3}$ and Proposition $\ref{rational}$.
\end{proof}
\begin{center}
{\bf Acknowledgements\\[1ex]}
\end{center}
The authors thank  Prof. B. Berndt and M. Beck  for communications
about reciprocity laws.

\end{document}